\documentclass[12pt,a4paper]{amsart}
\usepackage{amssymb, array, verbatim, amscd}
\usepackage{amsmath, amsfonts,amsthm,graphicx}

\usepackage{epsfig}
\usepackage{graphics}
\numberwithin{equation}{section}
\newcommand{\bc}{\begin{center}}
\newcommand{\ec}{\end{center}}
\newcommand{\bdm}{\begin{displaymath}}
\newcommand{\edm}{\end{displaymath}}
\newcommand{\beq}{\begin{equation}}
\newcommand{\eeq}{\end{equation}}
\newcommand{\bfl}{\begin{flushleft}}
\newcommand{\efl}{\end{flushleft}}
\newcommand{\bt}{\begin{tabbing}}
\newcommand{\et}{\end{tabbing}}

\newcommand{\Z}{{\mathbf Z}}

\newtheorem{theorem}{Theorem}[section]

\newtheorem{definition}[theorem]{Definition}
\newtheorem{definition2}[theorem]{Theorem and Definition}

\newcommand{\N}{{\mathbb N}}

\begin{document}
\title[Simulation of Probabilistic Sequential Systems ]{Simulation  of
Probabilistic Sequential Systems}
\author{Mar\'{\i}a A. Avi\~n\'o-Diaz}
\address{Cond. Parque San Antonio 1, Apart. 603\\
 Puerto Rico\\ Caguas, PR  00725 } \email{mavino2006@hotmail.com}
 \thanks{\today}


\thanks {These work  was partially supported by the Partnership
between M. D. Anderson Cancer Center, University of Texas and the
Medical Sciences Cancer Center, University of Puerto Rico, by  the
Program AABRE of Rio Piedras Campus,  and by the SCORE Program, NIH,
University of Puerto Rico en Cayey. }

 \maketitle

\begin{abstract}In this paper we introduce the idea of
probability in the definition of Sequential Dynamical Systems, thus
obtaining a new concept,  Probabilistic Sequential System. The
introduction of a probabilistic structure on Sequential Dynamical
Systems is an important and interesting problem.
 The notion of homomorphism of our new model, is a
natural extension of  homomorphism of sequential dynamical systems
introduced and  developed by Laubenbacher and Paregeis in several
papers. Our model, give the possibility to describe the dynamic of
the systems using Markov chains and all the advantage  of stochastic
theory. The notion of simulation  is introduced using the concept of
homomorphisms, as usual. Several examples of homomorphisms,
subsystems and simulations are given.
\end{abstract}
\section*{Introduction}\label{intro}
  \par
 Genetic Regulatory Networks had been modeled
using discrete and continuous mathematical models. An important
contribution to the simulation science is the theory of sequential
dynamical systems (SDS) \cite{B1, B2, B3,LP1,LP2}. In these paper,
the authors developed a new theory about the sequential aspect of
the entities  in a dynamical systems. In particular Laubenbacher
and Pareigis created  an elegant mathematical background of the SDS,
and with it solve several aspects of the theory  and applications.

Probabilistic Boolean Networks (PBN) had been  recently introduced
\cite{S1, S2, S3}, to model regulatory gene networks. The PBN are a
generalization of the widely used Boolean Network model (BN)
proposed by Kauffman (1969) \cite{K,K1}. While the PBN eliminate one
of the main limitations of the BN model, namely its inherent
determinism, they do not provide the framework for considering the
sequential behavior of the genes in the network, behavior observed
by the biologists.
\par

Here, we introduce the probabilistic structure in  SDS, using for
each vertices of the support graph a set of local functions, and
more than one schedule in the sequence of the local function
selected to form the update functions, obtaining a new concept:
probabilistic sequential dynamical system (PSS).  The notion of
simulation of a PSS is introduced in Section \ref{EHOM} using the
concept of homomorphism of PSS; and we prove
 that the  category of \textbf{SDS} is a full subcategory of the
of the category \textbf{ PSS}.  Several examples of homomorphisms,
subsystems and simulations are given.

 On the other hand, Deterministic Sequential Dynamical Systems have
 been studied for the last few years. The introduction of a probabilistic
 structure on Sequential Dynamical Systems is an important and interesting problem.
Our approach  take into account, a number of issues that have
already been recognized and solved for Sequential Dynamical Systems.
\section{Sequential Dynamical Systems and SDS-homomorphism}\label{SDS}
This section is an introduction with the definitions and results of
Sequential Dynamical System  introduced  by Laubenbacher and
Pareigis. Here we use   SDS  over a finite field. In this paper, we
denote the finite field $GF(p^{r})$ by $\mathcal{K}$, where $p$ is a
prime number.
\subsection{ Sequential Dynamical
Systems over finite fields}\label{defSDS} A Sequential Dynamical
System (\underline{SDS}) over a finite field $\mathcal {F}=(\Gamma,
  (f_i)_{i=1}^n, \alpha) $   consists of
\begin{enumerate}
\item [(1)]  A finite graph
$\Gamma=(V_\Gamma,E_\Gamma)$ with $V_\Gamma=\{1,\ldots ,n\}$
vertices, and a set of edges $E_\Gamma \subseteq V_\Gamma \times
V_\Gamma$.
 \item[(2)]  A family of  local functions $f_i:
{\mathcal{K}}^n \rightarrow {\mathcal{K}}^n $,  that is  \[f_i(x_1,
\ldots , x_n)= (x_1, \ldots ,x_{i-1}, \overline{f}, x_{i+1}, \ldots
, x_n)\] where $\overline{f}(x_1,\ldots , x_n)$ depends only of
those variables which are connected to $i$ in $\Gamma$.
\item [(3)] A
permutation $\alpha =( \begin{array}{lll} \alpha (1)& \ldots &
\alpha (n)
\end{array}) $ in the set of vertices $V_\Gamma$, called an update
schedule ( i.e. a bijective map $\alpha :V_\Gamma\rightarrow
V_\Gamma)$.
\end{enumerate}
\par The global update function of the \underline{SDS} is $f=f_{\alpha(1)}\circ\ldots\circ{f_{\alpha(n)}}$.
 The function $f$ defines the dynamical behavior of the \underline{SDS} and determines a finite directed graph
  with vertex set ${\mathcal{K}}^n$ and directed edges $(x,f(x))$, for all $x\in {\mathcal{K}}^n$, called the State
  Space of $\mathcal{F}$.

\subsection{Homomorphisms of \underline{SDS} }\label{homSDS}
The definition of homomorphism between two \underline{SDS} uses the
fact that the vertices $V_\Gamma=\{1,\ldots , n\}$ of an
\underline{SDS} and the  states $\mathcal{K}^n$ together with their
evaluation map $$\mathcal{K}^{n}\times V_{\Gamma} \ni (x,a) \mapsto
<x,a>:=x_a\in \mathcal{K}_a$$ form a contravariant setup, so that
morphisms between such structures should be defined contravariantly,
i.e. by a pair of certain maps $(\phi : \Gamma \rightarrow \Delta ,
h_\phi : K^{m}\rightarrow K^{n})$ or by a pair $(\phi : \Delta
\rightarrow \Gamma , h_\phi : K^{n}\rightarrow K^{m})$, with the
graph $\Delta$ having $m$ vertices. Here we use  a notation slightly
different the one using in  \cite{LP2}.

Let $\overline{F}=(\Gamma,(f_i:\mathcal{K}^n\rightarrow
\mathcal{K}^n),\alpha)$ and
$\overline{G}=(\Delta,(g_i:\mathcal{K}^m\rightarrow
\mathcal{K}^m),\beta)$ be two \underline{SDS}. Let $\phi:\Delta
\rightarrow \Gamma$ be a digraph morphism, and
\[(\widehat{\phi}_b :\overline{}\mathcal{K}_{\phi(b)}\rightarrow \mathcal{K}_b, \forall b\in \Delta),\]
 be a family of maps in the
category of \textbf{Set}. The map $h_\phi$ is an adjoint map,
because is defined as follows: consider the pairing
\[\mathcal{K}^n\times V_{\Gamma} \ni (x,a) \mapsto
<x,a>:=x_a\in \mathcal{K}_a;\] and similarly
\[\mathcal{K}^m\times V_\Delta \ni (x,b) \mapsto
<x,b>:=x_b\in \mathcal{K}_b.\] The induced  adjoint map is
$<h_\phi(x),b>:=\hat{\phi}_b(<x,\phi(b)>)=\hat{\phi}_b(x_{\phi(b)})$.
Then $\phi,$ and $(\widehat{\phi}_b)$ induce the adjoint map
$h_\phi:\mathcal{K}^{n}\rightarrow\mathcal{K}^{m}$ defined as
follows:
\begin{enumerate}\label{defh}
           \item [(I)]\hspace{.2in}
           $h_\phi(x_1,\ldots,x_n)=(\widehat{\phi}_{b_1}(x_{\phi_{(b_1)}}),\ldots,\widehat{\phi}_{b_m}(x_{\phi_{(b_m)}})).$
\end{enumerate}

Then $h_\phi$ is an homomorphism of \underline{SDS} if for a set of
orders $\tau _{\beta}$ associated to $\beta$ in the connected
components of $\Delta $, the map $h_\phi$ holds the following
condition and the commutative diagram.
\[(D1.a)\hspace{.3in}\left( \prod_{\beta_j\in \phi^{-1}(\alpha_i)}\prod_{j=l}^s (g_{\beta_j}) \right)\circ
h_\phi=h_\phi\circ f_{\alpha_i}\hspace{.4in}\begin{array}{lcl}
\mathcal{K}^n & \overset{f_{\alpha_i}}{-----\longrightarrow}& \mathcal{K}^n \\
\mid h_\phi&& \mid h_\phi\\
\downarrow  &  & \downarrow  \\
 \mathcal{K}^m &\underset{g_{\beta_l}\circ g_{\beta_{l+1}}\circ
\cdots \circ g_{\beta_s}} {-----\longrightarrow} & \mathcal{K}^m \cr
\end{array}\]
where $\prod_{\beta_j\in \phi^{-1}(\alpha_i)}\prod_{j=l}^s
(g_{\beta_j})=g_{\beta_l}\circ g_{\beta_{l+1}}\circ \cdots \circ
g_{\beta_s}$ and $\{\beta_l,\beta_{l+1},\ldots,\beta_s\}=\phi
^{-1}(\alpha_i)$.

 If $\phi^{-1}(\alpha_i)=\emptyset$, then
$Id_{\mathcal{K}^m} \circ h_\phi=h_\phi\circ f_{\alpha_i}$, and the
commutative diagram is the following.
\begin{enumerate}
\item [(D1.b)]\hspace{1.5in} $\begin{array}{lcl}
\mathcal{K}^n & \overset{f_{\alpha_i}}{-----\longrightarrow}& \mathcal{K}^n \\
\mid h_\phi&& \mid h_\phi\\
\downarrow  &  & \downarrow  \\
 \mathcal{K}^m &\underset{Id_{\mathcal{K}^m}} {-----\longrightarrow} & \mathcal{K}^m \cr
\end{array}$
\end{enumerate}

 For examples and properties see\cite{LP2}. It is clear
that the above diagrams implies the following one \begin{enumerate}
\item [(D2)]\hspace{1.7in}$\begin{array}{lcl}
\mathcal{K}^n & \overset{f=f_{\alpha_1}\circ \cdots \circ f_{\alpha_n}}{---\longrightarrow}& \mathcal{K}^n \\
\mid h_\phi&& \mid h_\phi\\
\downarrow  &  & \downarrow  \\
\mathcal{K}^m &\underset{g=g_{\beta_1}\circ \cdots \circ
g_{\beta_m}} {---\longrightarrow} & \mathcal{K}^m \cr
\end{array}$\end{enumerate}
\section{Probabilistic Sequential Dynamical Systems}
The following definition give us the possibility to have several
update functions acting in a sequential manner with assigned
probabilities. All these, permit us to study the dynamic of these
systems using Markov chains and other probability tools. We will use
the acronym PSS ( or \underline{SDS}) for plural as well as singular
instances.
\begin{definition}\label{PSS}
 A Probabilistic Sequential dynamical System (PSS) \[\mathcal
{D}_S=(\Gamma, \{F_i\}_{i=1}^{|\Gamma|=n}, (\alpha_j)_{j=1}^m,   C
)\] over $\mathcal{K}$ consists of:
\begin{enumerate}
\item[(1)] a finite graph
  $\Gamma=(V_\Gamma,E_\Gamma)$ with $n$ vertices;

  \item [(2)] a  set of local functions $F_i=\{f_{ik}:\mathcal{K}^n\rightarrow
\mathcal{K}^n|1\leq k \leq \ell(i)\}$ for each vertex $i$ of
$\Gamma$. (i. e. a bijection map $\sim : V_\Gamma\rightarrow \{F_i
|1\leq i\leq n\}$) (for definition of local function see
(\ref{defSDS}.2)).

\item [(3)] a family of $m$ permutations $\alpha_j  =(
\begin{array}{lll} \alpha_j (1)& \ldots & \alpha_j (n)
\end{array}) $ in the set of vertices $V_\Gamma $.

\item [(4)]and  a set $C=\{c_1, \ldots ,c_s\}$, of
selection probabilities.
\end{enumerate}
\end{definition}
We select one function  in  each set $F_i$, that is one for each
vertices $i$ of $\Gamma$, and a permutation $\alpha $, with the
order in which the vertices $i$ are selected (similarly SDS), so
there are $\underline{m}$ possibly different update functions
$f=f_{\alpha(1)k_{\alpha(1)}}\circ \ldots \circ
f_{\alpha(n)k_{\alpha(n)}}$, where $\underline{m}\leq m
!\times\ell(1)\times \ldots \times \ell (n)$. The probabilities are
assigned to the update functions, so there exists a subset
$S=\{f_1,\ldots ,f_s\}$ of update functions such that $c_k=p(f_k)$,
$1\leq k \leq s$.

The State Space of $\mathcal{D}_S$ is a digraph whose vertices are
the elements of $\mathcal{K}^{n}$ and there are an arrow going from
$x_1$ to $x_2$  if there exists   a functions $f$, such that
$x_2=f(x_1)$. For each one of the selected functions in $S$ we have
an SDS inside the PSS, so the state space of  the function $f$ is a
subdigraph of the state space of the PSS, so, the State Space of
$\mathcal{D}_S$  is a superposition of all inside  \underline{SDS}.
When we take the whole set of update functions generated by the
data, we will say that we have the \emph{full} PSS. We denote by
$\overline S$ the complement of set $S$ in the set of all  update
functions, and we will call the PSS $ {\mathcal {D}_{\overline S}}$
building by the same data but taking ${\overline S}$ as a set of
update functions, one complement of the PSS $\mathcal {D}_S$. All
the complements have the same set of function but they can use
different set of probabilities.

 The probability of the
arrow going from $x_1$ to $x_2$ for example, is the sum of the
probabilities of all functions $f$, such that $x_2=f(x_1)$. The PSS
represents a generalization of the \underline{SDS}: a
\underline{SDS} is a PSS for which every set of local functions has
one element,  and there is only one permutation in the family of
permutations. The update functions of the PSS have probabilities and
the state space of the PSS (or high level digraph) is described by a
transition matrix, and the dynamic is described by a Markov Chain.

\subsection{Example}\label{exam1}  Let
\[\mathcal{D}=(\Gamma;F_1,F_2,F_3;\alpha ^1, \alpha ^2;
F, (C(f_i))_{i=1}^8),\] be the following PSS:

(1) The graph: $\Gamma  \ \  \ \ \ \begin{array}{ccc} 1& \overset{\hspace{.2in}}{\overline{\hspace{.2in}}}& 3\\
& \diagdown &\mid  \\
 & & 2 \end{array}.$

(2)  If $\textbf{x}=(x_1,x_2,x_3)\in \{0,1\}^3$, then the sets of
local functions are the following:
\[\begin{array}{l}
F_1 =\{f_{11}(\textbf{x})=(1,x_2,x_3),f_{12}(\textbf{x})=(x_1+1,x_2,x_3))\}  \\
F_2 = \{f_{21}(\textbf{x})=(x_1, x_1x_2,x_3)\}   \\
F_3=\{f_{31}(\textbf{x})=(x_1,x_2,x_1x_2),f_{32}(\textbf{x})=(x_1,x_2,x_1x_2+x_3)\}
\end{array}\]

(3) The schedules or permutations are
$\alpha_1=\left(\begin{array}{lll} 3&2& 1
\end{array}\right); \alpha_2=\left(\begin{array}{lll} 1&2& 3
\end{array}\right) .$
  We obtain the following table of
 functions, and we select all of them for $\mathcal{D}$ because the
 probabilities given by $C$.
\[\begin{array}{ll}
f_1=f_{31}\circ f_{21} \circ f_{11} & f_2=f_{11}\circ f_{21} \circ f_{31} \\
f_3=f_{32}\circ f_{21} \circ f_{11} & f_4=f_{11}\circ f_{21} \circ
 f_{32}\\
f_5=f_{31}\circ f_{21} \circ f_{12} & f_6=f_{12}\circ f_{21} \circ
f_{31}\\
f_7=f_{32}\circ f_{21} \circ f_{12}& f_8=f_{12}\circ f_{21}\circ
f_{32} \cr
\end{array}.\]
  The  update functions are the following:
  \[\begin{array}{ll}
  f_1(\textbf{x})=(1,x_2,x_2)& f_2(\textbf{x})=(1,x_1x_2,x_1x_2) \\
 f_3(\textbf{x})=(1,x_2,x_2+x_3) & f_4(\textbf{x})=(1, x_1x_2,x_1x_2+x_3)\\
f_5(\textbf{x})=(1+x_1,(x_1+1)x_2,
 (x_1+1)x_2)  & f_6(\textbf{x})=(1+x_1,x_1x_2,x_1x_2)\\
f_7(\textbf{x})=(1+x_1,(x_1+1)x_2,(x_1+1)x_2+x_3)&
f_8(\textbf{x})=(1+x_1,x_1x_2,x_1x_2+x_3)
 \cr
\end{array}.\]
\begin{figure}
\includegraphics[height = 6in,width = 6in]{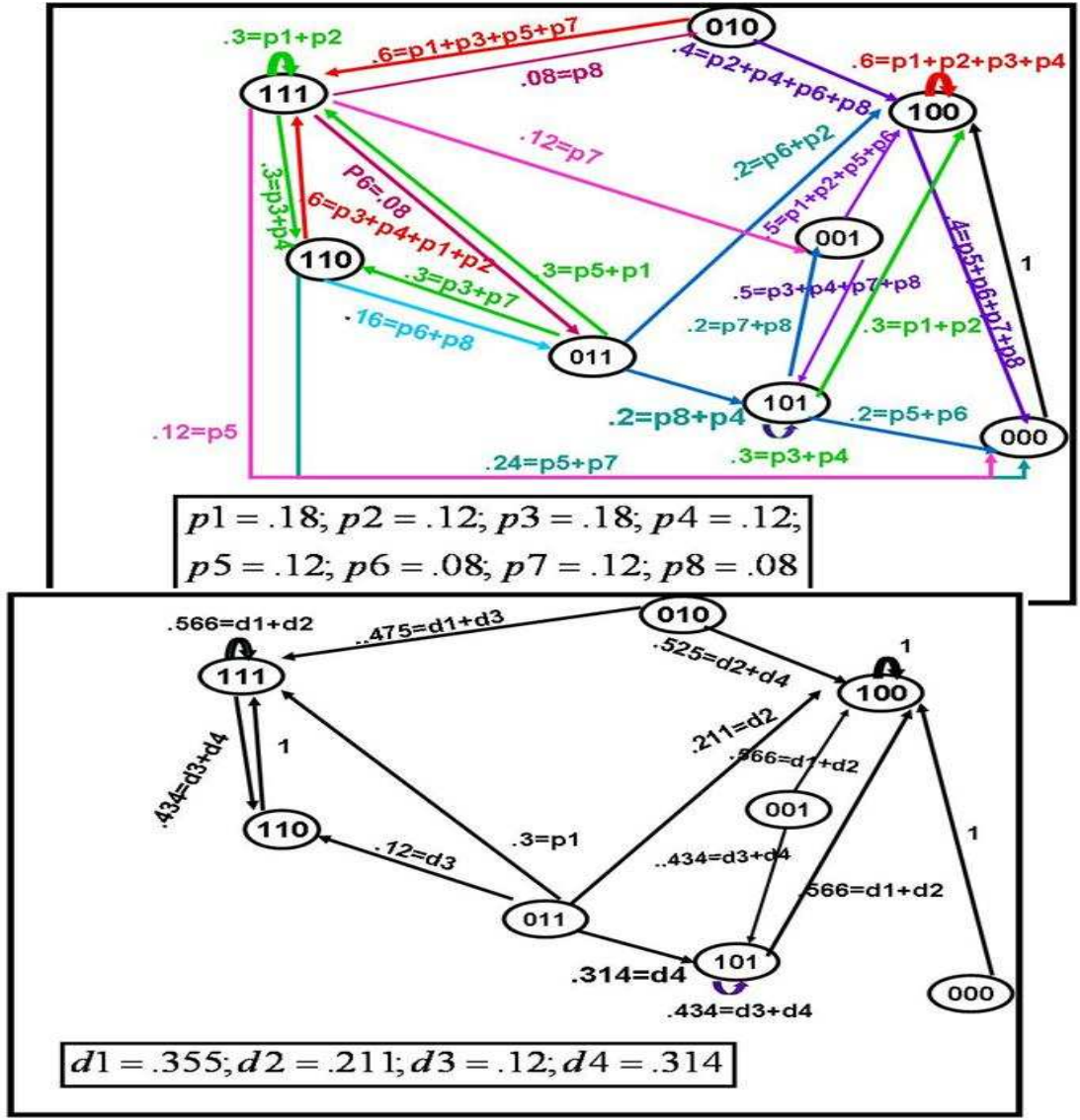}\caption{Above the State Space of $\mathcal{D}$, and below the State Space of $\mathcal{B}$ }
\end{figure}
(4) The  probabilities that we assign are: $p_1=C(f_1)=.18;
p_2=C(f_2)=.12; p_3=C(f_3)= .18; p_4=C(f_4)=.12; p_5=C(f_5)=.12;
p_6=C(f_6)=.08; p_7=c(f_7)=.12; p_8=c(f_8)=.08$. In order to study
the state space, it is convenient to determine
 the transition matrix of the system. The transition matrix of $\mathcal{D}$ is
$T$,
\[T=\left( \begin{array}{llllllll} 0&0&0&0&1&0&0&0\\
0&0&0&0&0.5&0.5&0&0\\
0&0&0&0&0.4&0&0&0.6\\
0&0&0&0&0.2&0.2&0.3&0.3\\
0.4&0&0&0&0.6&0&0&0\\
0.2&0.2&0&0&0.3&0.3&0&0\\
0.24&0&0&0.16&0&0&0&0.6\\
0.12&0.12&0.08&0.08&0&0&0.3&0.3
\end{array}\right).\]
\par
\subsubsection{Example}\label{MA}
We notice that there are several PSS that we can construct with the
same initial data of functions and permutations, but with different
set of  probabilities, that is,  sub-PSS of $\mathcal{D}$ see
\ref{homPSS}. For example if $S=\{f_1, f_2, f_3,f_4\}$ and $C=\{
d(f_1)=.355, d(f_2)=.211, d(f_3)=.12,d(f_4)=.314)\}$, then
\[\mathcal{B}=(\Gamma;F_1,F_2,F_3;\alpha _1, \alpha _2; S; C ).\]
\section{Homomorphisms of Probabilistic Sequential Systems}
One objective in this section is to show that: homomorphism of PSS
is a natural extension of homomorphism of \underline{SDS}. We
remark, in the introductory  section \ref{SDS}, the definition  of
\underline{SDS}. In the definition of a homomorphism of PSS we
establish conditions to connect the support graphs, the State Space
and the assigned probabilities.

\subsection{Definition of homomorphisms of PSS}\label{homPSS}
Let $\mathcal{D}_1=(\Gamma,(F_i)_{i=1}^{|\Gamma|=n},(\alpha ^j)_{j},
S_1,C_1)$ and $\mathcal{D}_2
=(\Delta,(G_i)_{i=1}^{|\Delta|=m},(\beta ^k)_{k},S_2,C_2)$ be two
PSS over a finite field $\mathcal{K}$.

A homomorphism  from $\mathcal{D}_1$ to $\mathcal{D}_2$ is a pair of
functions $(\phi, h_\phi)$ where:
\begin{enumerate}
\item [(\ref{homPSS}.1)]  $\phi:\Delta \rightarrow \Gamma $ is a
graph morphism, and
\[(\widehat{\phi}_{b}:\mathcal{K}_{\phi(b)}\rightarrow \mathcal{K}_{b}\forall b\in \Delta ),\] is a family of maps in the
category of \textbf{Set}. They induce the adjoint function $h_\phi$,
see (\ref{defh}) (I).
 \item [(\ref{homPSS}.2)]  The induced adjoint map $h_\phi: \mathcal{K}^{n}\rightarrow \mathcal{K}^{m}$ is
a map such that for all update function $f$ in $S_1$ there exists an
update function $g\in S_2$ such that $(\phi, h_\phi)$ is an
\underline{SDS}-homomorphism from
$(\Gamma,(f:\mathcal{K}^n\rightarrow \mathcal{K}^n),\alpha^j)$ to
$(\Delta,(g:\mathcal{K}^m\rightarrow \mathcal{K}^m),\beta ^k)$. That
is, the  diagrams (D1), and (D2) commute for all  $f$ and the
selected $g$.
 \item [(\ref{homPSS}.3)](The $\epsilon$-condition)
 For a fixed real
number $0\le \epsilon<1$,  the map $h_\phi$ satisfies the following:
\begin{equation*}\label{epsilon}
 max_{u,v}|c_{f_i}(u,v)-d_{g_j}(h_\phi(u),h_\phi(v))|\le \epsilon
 \end{equation*}
 for all  $f$ in $S_1$, and its selected
$g$ in $S_2$ by  condition (\ref{homPSS}.2), and $u,\ v\in
\mathcal{K}^n$,
\end{enumerate}
We will say that an homomorphism $(\phi, h_\phi)$ from
$\mathcal{D}_1$ to $\mathcal{D}_2$ is an \textbf{isomorphism} if
$\phi$ and $h_\phi$ are bijective functions, and
$d(h_\phi(u),h_\phi(g(u))= c(u,f(u))$ for all $u$, in
$\mathcal{K}^n$, and all  $f$. We denote it by $\mathcal{D}_1\cong
\mathcal{D}_2$.

\subsection{Theorem}\label{C3}\emph{In the definition of homomorphism the
condition (\ref{homPSS}.3) is a consequence of the condition
(\ref{homPSS}.1), and (\ref{homPSS}.2). }
\begin{proof}
Suppose $(\phi,h_\phi)$ satisfies (\ref{homPSS}.1),
(\ref{homPSS}.2), and
$max_{u,v}|c_{f_i}(u,v)-d_{g_j}(h_\phi(u),h_\phi(v))|\ge 1$, then
for some $(u,v)$ we have one of the following cases
\begin{itemize}
\item [(Case 1)] $c_{f_i}(u,v)=1$ and
$d_{g_j}(h_\phi(u),h_\phi(v))=0$,  it is impossible by
(\ref{homPSS}.2).
\item[(Case 2)] $d_{g_j}(h_\phi(u),h_\phi(v))=1$, and
$c_{f_i}(u,v)=0$, it is impossible because there exists at least
other element $w\in \mathcal{K}^{n}$ such that $c_{f_i}(u,w)\ne 0$,
and so on $d_{g_j}(h_\phi(u),h_\phi(w))\ne 0$. Because the sum of
probabilities of all arrow going up of $h_\phi(u)$ is equal $1$,
then $d_{g_j}(h_\phi(u),h_\phi(v))<1$.
\end{itemize}
Therefore the third   condition (\ref{homPSS}.3) holds, and always
$\epsilon$ exists.
\end{proof}
 \subsection{Theorem}\label{MT}
\emph{If $(\phi, h_\phi):\mathcal{D}_1\longrightarrow\mathcal{D}_2$
is a PSS-homomorphism, and $T_i$ denote the transition matrix of
$\mathcal{D}_i$, and the entry $(u,v)$ of $T_1$ is $p(u,v)$ then:
\[\lim
_{m\rightarrow\infty}|({T_1})^{m}_{u,v}-({T_2})^{m}_{h_\phi(u),h_\phi(v)}|=0,\]
 for all possible  $u$ and $v$ in $\mathcal{K}^n$.}
\begin{proof}
The $\epsilon$-condition (\ref{homPSS}.3) asserts that
 for a fixed real
number $0\le \epsilon<1$,  the map $h_\phi$ satisfies the following:
\begin{equation*}\label{epsilon}
 max_{u,v}|c_{f_i}(u,v)-d_{g_j}(\phi(u),\phi(v))|\le \epsilon
 \end{equation*}
 for all  $f$ in $S_1$, and its selected
$g$ in $S_2$ by  condition (\ref{homPSS}.2), and $u,\ v\in
\mathcal{K}^n$.

If we  have a function $f_i$ going from $u$ to $v=f_i(u)$  in
$\mathcal{K}^{n}$, then there exists a function $g_j$ going from
$h_\phi(u)$ to $h_\phi(v)$, so  $g_j(h_\phi(u))=h_\phi(f_i(u))$.
Then for $m= 2$,
 and by the Chapman-Kolmogorov equation \cite{Ste,KO}, we
have the following:
\[|c_{f_i^2}(u,f_i^2(u))-d_{g_j^2}(h_\phi(u),g_j^2(h_\phi(u)))|=\]
\[|c_{f_i}(u,f_i(u))c_{f_i}(f_i(u),f_i^2(u))-d_{g_j}(h_\phi(u),g_j(h_\phi(u)))d_{g_j}(g_j(h_\phi(u)),g_j^2(h_\phi(u)))|=\]
\[|c_{f_i}(u,f_i(u))c_{f_i}(f_i(u),f_i^2(u))-d_{g_j}(h_\phi(u),h_\phi(f_i(u)))d_{g_j}(h_\phi(f_i(u)),h_\phi(f_i^2(u)))|\leq\]
\[\leq
|c_{f_i}(f_i(u),f_i^2(u))|\epsilon+|d_{g_j}(h_\phi(u),h_\phi(f_i(u)))|\epsilon\leq
2\epsilon.\] By condition (2) in definition of homomorphism. Then we
proved that
$|c_{f_i^2}(u,f_i^2(u))-d_{g_j^2}(\phi(u),g_j^2(\phi(u)))|\leq 2
\epsilon$.
 Using  mathematical induction  over $m$, we can conclude that
\begin{itemize}
  \item [(\ref{MT}.1)] \hspace{.6in}$\max|(c_{f_i^m}(u,f_i^m(u))-d_{g_j^m}(\phi(u),g_j^m(\phi(u)))|\leq
m\epsilon $
\end{itemize}
 Let $S_\phi$ be the set of function in $S_2$ associated to a function in $S_1$.
 If $u,\ v \in \mathcal{K}^n$, and we denote by
$p_2(u,v)=\sum_{f^{2}_i}c_{f^{2}_i}(u,v)$ and
$p_2(h_\phi(u),h_\phi(v))=\sum_{g^{2}_j}d_{g^{2}_j}(h_\phi(u),h_\phi(v))$,
then condition (\ref{homPSS}.3) implies that
$|p_2(u,v)-p_2(h_\phi(u),h_\phi(v))|\leq 2k \epsilon + \delta^{2}$,
where $k$ is the maximum number of functions $f^{2}_i$ going from
one state to another in the state space $\mathcal{K}^n$. Denoting
$\delta^{2}=\sum_{g^{2}_j \not \in S_\phi}
d_{g^{2}_j}(h_\phi(u),h_\phi(v)) $we know that $\delta^{2} < 1$. So,
if $T_i$ denote the transition matrix of $\mathcal{D}_i$, and the
entry $(u,v)$ of $(T_1)^{2}$ is $p_2(u,v)$ then the condition
(\ref{homPSS}.3) implies that:
\begin{itemize}
  \item [(\ref{MT}.2)] \hspace{1in}$\max_{u,v}|({T_1})^{2}_{u,v}-({T_2})^{2}_{h_\phi(u),h_\phi(v)}|\leq
2k\epsilon + \delta ^{2}$\end{itemize}  for all possible  $u$ and
$v$ in $\mathcal{K}^n$.

Using equations (\ref{MT}.1), and  (\ref{MT}.2), and denoting
$\delta^{m}=\sum_{g_j \not \in S_\phi}
[d_{g_j}(h_\phi(u),h_\phi(v))]^{m} $,  we conclude that
\begin{itemize}
  \item [(\ref{MT}.3)] \hspace{1in}$
\max_{u,v}|({T_1})^{m}_{u,v}-({T_2})^{m}_{h_\phi(u),h_\phi(v)}|\leq
mk\epsilon + \delta ^{m}$ \end{itemize}   for all possible $u$ and
$v$ in $\mathcal{K}^n$.

So, for all real number $0<\epsilon ' < 1$ there exists  $m\in\N$,
such that,
\[|({T_1})^{m'}_{u,v}-({T_2})^{m'}_{h_\phi(u),h_\phi(v)}|< \epsilon ',\]
 for all natural number  $m'>m$, and for all
possible $u, \ v \in \mathcal{K}^n$.

 In fact, using notation of equation (\ref{MT}.2), we have $\epsilon ^{m'}\ll \epsilon
 ^{m}$, and this implies
 $(m'k)\epsilon ^{m'}<(mk) \epsilon ^{m}$. Similarly
  $\delta ^{m'}< \delta ^{m}$, So, selecting $m$ such that
   $(mk)\epsilon ^{m}+ \delta ^{m}< \epsilon '$, we obtain
\[|({T_1})^{m'}_{u,v}-({T_2})^{m'}_{h_\phi(u),h_\phi(v)}|\leq (m'k) \epsilon ^{m'}+ \delta ^{m'}< (mk)\epsilon^{m}+ \delta ^{m}< \epsilon ',\]
where $k$ is the maximum number of functions going from one state to
another in the state space of the power $m'$ of the functions .
Therefore
\[\lim
_{m\rightarrow\infty}|({T_1})^{m}_{u,v}-({T_2})^{m}_{h_\phi(u),h_\phi(v)}|=0,\]
 for all possible  $u$ and $v$ in $\mathcal{K}^n$, and
 the theorem holds
\end{proof}

\subsection{Definition of $\epsilon$-equivalent}\label{equiv} If $\phi$ and $h_\phi$ are bijective functions,
and the condition (\ref{homPSS}.3) holds, but the probabilities are
not equal, we will say that $\mathcal{D}_1$, and  $\mathcal{D}_2$
are $\epsilon$-equivalent, and we write $\mathcal{D}_1 \simeq
_\epsilon \mathcal{D}_2$. So, $\mathcal{D}_1$, and  $\mathcal{D}_2$
are $\epsilon$-equivalent if there exist  $(h_\phi,\phi)$, and
$((h_\phi^{-1},\phi^{-1})$, such that for all $f\in S_1$ and $g\in
S_\phi$, we have $f=h_\phi^{-1}\circ g \circ h_\phi$.

\subsection{Image of $h_\phi$}\label{prob} Consider $\underline{d}(g)$  defined as follows: if
$S_\phi=\{g_1,\ldots, g_s\}$ is  the set of  functions in $S_2$
selected in (\ref{homPSS}.2) for the map $h_\phi$  then the new
probability of $g_i$ is
\[
 \underline{d}(g_i)=\frac{d(g_i)}{\sum_{g_i\in S_\phi}d(g_i)}.
\]
 The set of functions  $S_\phi$ in $\mathcal{G}$, together with the new
probabilities defined above, form  a new PSS, that we will call
Image of $h_\phi=Im (h_\phi)$. So, the graph $\Delta$, the update
functions $S_\phi$ determine de local functions associate to each
vertex of $\Delta$, similarly the permutations using by these
functions, and finally the new probabilities assigned.
\subsection{Sub Probabilistic Sequential System}\label{subPSS}
We will say an injective monomorphism is a PSS-homomorphism  such
that $\phi$ is surjective and the set of functions $\hat{\phi}_i$,
for all $i$ are injective functions, and so $h_\phi$, \cite{LP2}.
Therefore, we will say that a PSS $\mathcal{G}_X$ is sub
Probabilistic Sequential System of $\mathcal{F}_S$ if there exists
an injective monomorphism from $\mathcal{G}_X$ to $\mathcal{F}_S$.
\subsection{Definitions}\label{DEFH}\hspace{2.5in}\par

(\ref{DEFH}.1) Let $\mathcal{F}=(\Gamma,(F_i)_{i=1}^n,(\alpha
^j)_{j\in{J}},C)$ be a PSS. The pair of functions
$\mathcal{I}=(id_\Gamma,id_{\mathcal{K}^n})$ is the \emph{identity
homomorphism}, and it is an example of an isomorphism.

(\ref{DEFH}.2) An homomorphism $(\phi_1, h_{\phi_1})$ of PSS is an
\emph{injective monomorphism} if $\phi$ is surjective and
$h_{\phi_1}$ is injective, for example see (\ref{EHOM}.3). Similarly
we will say that an homomorphism is a \emph{surjective epimorphism}
if $\phi$ is injective and $h_{\phi_1}$ is surjective,  for a
complete description of the properties of this class of monomorphism
and epimorphism see Section 7 in \cite{LP2}.

\subsection{Simulation of  PSS}
 We consider that the PSS $\mathcal{G}$ is simulated by $\mathcal{F}$ if
there exist a injective monomorphisms $(\phi_1, h_{\phi_1}):
F\rightarrow G$ or a surjective epimorphism $(\phi_2,
h_{\phi_2}):\mathcal{G}\rightarrow \mathcal{F}$.

 \section{Examples of Simulation and  Homomorphisms of PSS. }\label{EHOM}
In this section we give several examples of PSS-homomorphism, and
simulation. In the second example we show how the condition
(\ref{homPSS}.2) is verified under the supposition that a function
$\phi$ is defined. So, we have two examples in (\ref{EHOM}.2), one
with $\phi$ the natural inclusion, and the second  with $\phi$ the
only possibility of a surjective map. In the last example we have a
complete description of two PSS, where we have only one permutation
and two or less functions for each vertices in the graph. In
particular this homomorphism is an injective monomorphism, so is an
example  of simulation too.

 \textbf{(\ref{EHOM}.1)} For the PSS, in the examples \ref{exam1}, and \ref{MA} we now
define the natural inclusion $(Id_\Gamma, \iota):
\mathcal{B}\rightarrow \mathcal{D}$. It is clear that the inclusion
satisfies the two first condition to be a homomorphism, and the
third one is a simple consequence of the Theorem \ref{C3}. In fact
\[T_{\mathcal{B}}\left( \begin{array}{llllllll} 0&0&0&0&1&0&0&0\\
0&0&0&0&0.566&0.434&0&0\\
0&0&0&0&0.525&0&0&0.475\\
0&0&0&0&0.211&0.314&0.12&0.355\\
0&0&0&0&1&0&0&0\\
0&0&0&0&0.566&0.434&0&0\\
0&0&0&0&0&0&0&1\\
0&0&0&0.08&0&0&0.434&0.566
\end{array}\right)\]
and $|(T_{\mathcal{B}})_{u,v}-(T)_{\iota(u),\iota(v)}|<0.4$. Here
${\mathcal{D}}$ is simulated by ${\mathcal{B}}$

 \textbf{(\ref{EHOM}.2)} Consider the two
graphs below
\[\Gamma \ \  \begin{array}{ccc} 2&\overset{\hspace{.2in}}{\overline{\hspace{.2in}}} & 3\\
\mid &   &\  \\
1&  \overset{\hspace{.2in}}{\overline{\hspace{.2in}}} & 4
\end{array} \hbox{   and       } \Delta \  \  \begin{array}{ccc} 2&\overset{\hspace{.2in}}{\overline{\hspace{.2in}}} & 3   \\
\mid &   & \\
1& &
\end{array}\]
 Suppose that the functions
associated to the vertices are the  families $\{f_1,f_2,f_3,f_4\}$,
 for $\Gamma $ and $\{g_1,g_2,g_3\}$ for $\Delta $. The
permutations are $\alpha_1=(4\ 3\ 2\ 1)$, $\alpha_2=(4\ 1\ 3\ 2)$
and $\beta_1=(3\ 2\ 1)$, $ \beta_2=(1\ 3\ 2)$, so, $S=\{f=f_4\circ
f_3 \circ f_2 \circ f_1; \underline{f}=f_4\circ f_1 \circ f_3 \circ
f_2\}$, and $X=\{g= g_3 \circ g_2 \circ g_1; \underline{g}= g_1
\circ g_3 \circ g_2\}$. Then, we  have constructed two PSS, each one
with two permutations and only one function associated to each
vertex in the graph; denoted by:
\[\mathcal{D}_S=(\Gamma;f_1,f_2,f_3,f_4; \alpha^1, \alpha^2;S;C ) \ and \
\mathcal{B}_X=(\Delta;g_1,g_2,g_3;\beta^1, \beta^2;X;D ).\]

\emph{Case (a)} We assume that there exists a homomorphism $(\phi ,
h_{\phi})$ from $\mathcal{D}_S$ to
 $\mathcal{B}_X$, with the graph morphism  $\phi: \Delta
\rightarrow \Gamma $ is given by $\phi(1)=1,\ \phi(2)=2, \
\phi(3)=3$. Suppose the functions \[(\widehat{\phi}_b
:\mathcal{K}_{\phi(b)}\rightarrow \mathcal{K}_b, \forall b\in
\Delta),\] are giving, and the adjoint function
\[h_{\phi}:{\Z_p}^4\rightarrow {\Z_p}^3, \  h_\phi (x_1,x_2,x_3,x_4)=(\hat{\phi}_1(x_1),\hat{\phi}_2(x_2),\hat{\phi}_3(x_3))\] is defined too.
If $(h_\phi,\phi)$ is an homomorphism, which satisfies the
definition (\ref{homPSS}), then the following diagrams commute:
\[\begin{array}{ccccccccc}
{\Z_p}^4 &\overset{f_4}{\longrightarrow}&{\Z_p}^4
&\overset{f_3}{\longrightarrow} & {\Z_p}^4
&\overset{f_2}{\longrightarrow} & {\Z_p}^4
&\overset{f_1}{\longrightarrow} & {\Z_p}^4  \\
 h_{\phi}\downarrow  & &h_{\phi} \downarrow &  & h_{\phi}\downarrow & & h_{\phi}\downarrow &
&h_{\phi}\downarrow \\
{\Z_p}^3 &\overset{Id}{\longrightarrow} &{\Z_p}^3
&\overset{g_3}{\longrightarrow} & {\Z_p}^3
&\overset{g_2}{\longrightarrow} & {\Z_p}^3
&\overset{g_1}{\longrightarrow} & {\Z_p}^3
\end{array}, \  \begin{array}{ccc}
 {\Z_p}^4 &\overset{f}{\longrightarrow}&{\Z_p}^4\\
  h_{\phi}\downarrow  & &h_{\phi} \downarrow \\
  {\Z_p}^3 &\overset{g}{\longrightarrow} & {\Z_p}^3 \cr
  \end{array} \]
\[\begin{array}{ccccccccc}
{\Z_p}^4 &\overset{f_4}{\longrightarrow}&{\Z_p}^4
&\overset{f_1}{\longrightarrow} & {\Z_p}^4
&\overset{f_3}{\longrightarrow} & {\Z_p}^4
&\overset{f_2}{\longrightarrow} & {\Z_p}^4  \\
 h_{\phi}\downarrow  & &h_{\phi} \downarrow &  & h_{\phi}\downarrow & & h_{\phi}\downarrow &
&h_{\phi}\downarrow \\
{\Z_p}^3 &\overset{Id}{\longrightarrow} &{\Z_p}^3
&\overset{g_1}{\longrightarrow} & {\Z_p}^3
&\overset{g_3}{\longrightarrow} & {\Z_p}^3
&\overset{g_2}{\longrightarrow} & {\Z_p}^3
\end{array}, \  \begin{array}{ccc}
 {\Z_p}^4 &\overset{\underline{f}}{\longrightarrow}&{\Z_p}^4\\
  h_{\phi}\downarrow  & &h_{\phi} \downarrow \\
  {\Z_p}^3 &\overset{\underline{g}}{\longrightarrow} & {\Z_p}^3 \cr
  \end{array}.\]

\emph{Case (b)} Consider now the map $\phi: \Gamma \rightarrow
\Delta$, defined by $\phi(1)=1$, $\phi(2)=2$, $\phi(3)=3$, and
$\phi(4)=1$. If there exists an homomorphism  $(\phi , h_{\phi}):
\mathcal{B}_X\rightarrow \mathcal{D}_S$ that satisfies
(\ref{homPSS}.2, then

\[\begin{array}{ccccccc}
 {\Z_p}^3 &\overset{g_3}{\longrightarrow} & {\Z_p}^3 &\overset{g_2}{\longrightarrow}& {\Z_p}^3
&\overset{g_1}{\longrightarrow} & {\Z_2}^3  \\
  h_{\phi}\downarrow & & h_{\phi}\downarrow& & h_{\phi}\downarrow &
&h_{\phi}\downarrow \\
 {\Z_p}^4 &\overset{f_4\circ f_3}{\longrightarrow} & {\Z_p}^4 &\overset{f_2}{\longrightarrow}
& {\Z_p}^4 &\overset{f_1}{\longrightarrow} & {\Z_p}^3
\end{array}, \  \begin{array}{ccc}
 {\Z_p}^3 &\overset{g}{\longrightarrow}&{\Z_p}^3\\
  h_{\phi}\downarrow  & &h_{\phi} \downarrow \\
  {\Z_p}^4 &\overset{f}{\longrightarrow} & {\Z_p}^4 \cr
  \end{array} \]
\[\begin{array}{ccccccc}
 {\Z_p}^3 &\overset{g_1}{\longrightarrow} & {\Z_p}^3 &\overset{g_3}{\longrightarrow}& {\Z_p}^3
&\overset{g_2}{\longrightarrow} & {\Z_2}^3  \\
  h_{\phi}\downarrow & & h_{\phi}\downarrow& & h_{\phi}\downarrow &
&h_{\phi}\downarrow \\
 {\Z_p}^4 &\overset{f_4\circ f_1}{\longrightarrow} & {\Z_p}^4 &\overset{f_3}{\longrightarrow}
& {\Z_p}^4 &\overset{f_2}{\longrightarrow} & {\Z_p}^3
\end{array}, \  \begin{array}{ccc}
 {\Z_p}^3 &\overset{\underline{g}}{\longrightarrow}&{\Z_p}^3\\
  h_{\phi}\downarrow  & &h_{\phi} \downarrow \\
  {\Z_p}^4 &\overset{\underline{f}}{\longrightarrow} & {\Z_p}^4 \cr
  \end{array} \]

\textbf{(\ref{EHOM}.3)} We now construct a PSS-homomorphism from
$\mathcal {F}_S=(\Gamma ,(F_i )_3,\beta,S ,C )$ to $\mathcal
{G}_X=(\Delta ,(G_i )_4,\alpha, X ,D)$, with the property that
$\phi$
 is surjective and the functions $\phi _i$ are injective, that we  call a injective monomorphism.
   The PSS $\mathcal {F}_S$ has a support graph $\Gamma$ with three vertices, and the PSS
$\mathcal {G}_X$ has a support graph $\Delta $ with four vertices
\[ \Gamma \  \ \begin{array}{ccc}
&& 3\\
&  & \mid \\
 1 & \overset{\hspace{.2in}}{\overline{\hspace{.2in}}} & 2
\end{array} \ \Delta \ \ \ \ \
\begin{array}{ccc}
 2&\overset{\hspace{.2in}}{\overline{\hspace{.2in}}} & 4\\
& \diagdown  & \mid \\
 1 & \overset{\hspace{.2in}}{\overline{\hspace{.2in}}} & 3
\end{array},  \hspace{.2in}   \    \]
The homomorphism $(h_\phi, \phi ):\mathcal {F}_S \rightarrow
\mathcal {G}_X$, has the contravariant  graph morphism $\phi :
\Delta\rightarrow \Gamma$, defined by the arrows of graphs, as
follows $\phi(1)=1$, $\phi(2)=\phi(3)=2,$ and $\phi(4)=3 $. The
family of functions $\hat{\phi}_i:\mathcal{K}_{\phi (i)}\rightarrow
 \mathcal{K}_{(i)}$,  $\hat{\phi}_1(x_1)=x_1$; $\hat{\phi}_2(x_2)=x_2;\
 \hat{\phi}_3(x_2)=x_2;\  \hat{\phi}_4(x_4)=x_4.$ The adjoint function is
 \[h_\phi:{\Z_2}^3\rightarrow {\Z_2}^4,\ \textrm{with} \  \  h_\phi(x_1,x_2,x_3)=(\hat{\phi}_1(x_1),\hat{\phi}_2(x_2),\hat{\phi}_3(x_2),\hat{\phi}_4(x_4))=(x_1,x_2,x_2,x_3).\]
 The first condition in the definition \ref{homPSS} holds.

The  PSS  $\mathcal{F}_S=(\Gamma ;(F_i)_3 ;\beta ;S ; C  ),$ with
data of functions $ F_1=\{ f_{11}=Id,
f_{12}(x_1,x_2,x_3)=(1,x_2,x_3)\},
 F_2=\{f_{21}(x_1,x_2,x_3)=(x_1,x_2,x_3)\},$ and $
F_3=\{f_{31}(x_1,x_2,x_3)=(x_1,x_2,x_2\overline{x_3})\},$ one
permutation or schedule  $\beta =(\ 1 \ 2 \ 3\ );$ and probabilities
$C=\{c(\underline{f})=.5168, c(\underline{g})=.4832\} $, so $S
=\{\underline{f},\underline{g}\}$; where the two update functions
are
\[\underline{f}= f_{11}\circ f_{21}\circ f_{31},\
\underline{f}(x_1,x_2,x_3)= (x_1,x_2,x_2\overline{x_3});\]
\[\textrm{and} \ \ \ \underline{g}= f_{12}\circ f_{21}\circ f_{31}, \
\underline{g}(x_1,x_2,x_3)= (1,x_2,x_2\overline{x_3}).\]

The PSS  $\mathcal {G}_X=(\Delta ;(G_i)_4;\alpha; X; D )$ is a PSS,
with the following data: the families of functions: $G_1=\{g_{11},
g_{12}\}$; $G_2=\{g_{21}, g_{22}\}$, $G_3=\{g_{31},g_{32}\}$; and
$G_4=\{g_4\}$.
\[\begin{array}{ll}
g_{11}(x_1,x_2,x_3,x_4)&= (1,x_2,x_3,x_4) \\
g_{21}(x_1,x_2,x_3,x_4)&= (x_1,1,x_3,x_4)  \\
g_{31}(x_1,x_2,x_3,x_4)&= (x_1,x_2,x_1x_2,x_4) \\
g_{41}(x_1,x_2,x_3,x_4)&=(x_1,x_2,x_3,x_2\overline{x_4})\\
\end{array}\ \begin{array} {ll}
g_{12}= Id=g_{22} &  \\
g_{32}(x_1,x_2,x_3,x_4)&= (x_1,x_2,x_2,x_4)\\
\end{array}.\]
One schedule   $\alpha=\left( \begin{array}{llll} 1&2& 3&
4\end{array}\right)$,  the eight possible  update functions,  and
its probabilities $D=\{d(g)=0,d(f)=0,
d(\tilde{f})=0,d(\tilde{g})=0,d(\hat{f})= .00252,d(\hat{g})=.08321,
d(\check{f})=.51428,d(\check{g})=.39999\}$  whose determine
$X=\{\hat{f},\hat{g},\check{f},\check{g}\}$ probabilities, are the
following:
\[\begin{array}{ll}
1)g(x_1,x_2,x_3,x_4)=(g_{11} \circ g_{21} \circ g_{31 }\circ g_4 )(x_1,x_2,x_3,x_4)=(1,1,x_1x_2,x_2\overline{x_4}),& \\
2)f(x_1,x_2,x_3,x_4)= (g_{12} \circ g_{21} \circ g_{32}\circ g_4)(x_1,x_2,x_3,x_4)=(x_1,1,x_2,x_2\overline{x_4}),& \\
3)\tilde{f}(x_1,x_2,x_3,x_4)=(g_{12} \circ g_{21} \circ g_{31}\circ g_4)(x_1,x_2,x_3,x_4)=(x_1,1,x_1x_2,x_2\overline{x_4}),&\\
4) \tilde{g}(x_1,x_2,x_3,x_4)=(g_{11}\circ g_{21}\circ g_{32}\circ g_4)(x_1,x_2,x_3,x_4)=( 1,1,x_2,x_2\overline{x_4}),&\\
5)\hat{f}(x_1,x_2,x_3,x_4)=(g_{12} \circ g_{22} \circ g_{31}\circ g_4)(x_1,x_2,x_3,x_4)=(x_1,x_2,x_1x_2,x_2\overline{x_4}),&\\
6)\hat{g}(x_1,x_2,x_3,x_4)=(g_{11} \circ g_{22} \circ g_{31} \circ g_4)(x_1,x_2,x_3,x_4)=(1,x_2,x_1x_2,x_2\overline{x_4}),& \\
7) \check{f}(x_1,x_2,x_3,x_4)=(g_{12} \circ g_{22} \circ g_{32}\circ g_4)(x_1,x_2,x_3,x_4)= (x_1,x_2,x_2,x_2\overline{x_4}),& \\
8)\check{g}(x_1,x_2,x_3,x_4)=(g_{11}\circ g_{21}\circ g_{32}\circ g_4)(x_1,x_2,x_3,x_4)= ( 1,x_2,x_2,x_2\overline{x_4}),&\\
\end{array}.\]
  We claim $(\phi,h_\phi):\mathcal{F}_S\rightarrow\mathcal{G}_X $ is a homomorphism. We will prove that the following diagrams commute.
\[\begin{array}{lcl}
{\Z_2}^3 &\overset{\underline{f}} {\longrightarrow}& {\Z_2}^3 \\
h_\phi \downarrow &  & \downarrow h_\phi \\
 {\Z_2}^4 & \overset{\check{f}}{\longrightarrow} & {\Z_2}^4 \cr
\end{array}, \ \textrm{and} \  \ \begin{array}{lcl}
{\Z_2}^3 & \overset{\underline{g}}{\longrightarrow}& {\Z_2}^3 \\
h_\phi \downarrow &  & \downarrow h_\phi \\
 {\Z_2}^4 &\overset{\check{g}}{\longrightarrow} & {\Z_2}^4 \cr
\end{array},\]
In fact, \[(h_\phi \circ \underline{f})(x_1,
x_2,x_3)=h_\phi(x_1,x_2,x_2\overline{x_3})
=(x_1,x_2,x_2,x_2\overline{x_3})=(\check{f}\circ h_\phi)(x_1,
x_2,x_3).
\] On the other hand,
\[(h_\phi \circ
\underline{g})(x_1,x_2,x_3)=h(1,x_2,x_2\overline{x_3})=(1,x_2,x_2,x_2\overline{x_3})=(\check{g}\circ
h_\phi)(x_1,x_2,x_3).\]
\begin{figure}
\includegraphics[height = 3in,width = 6in]{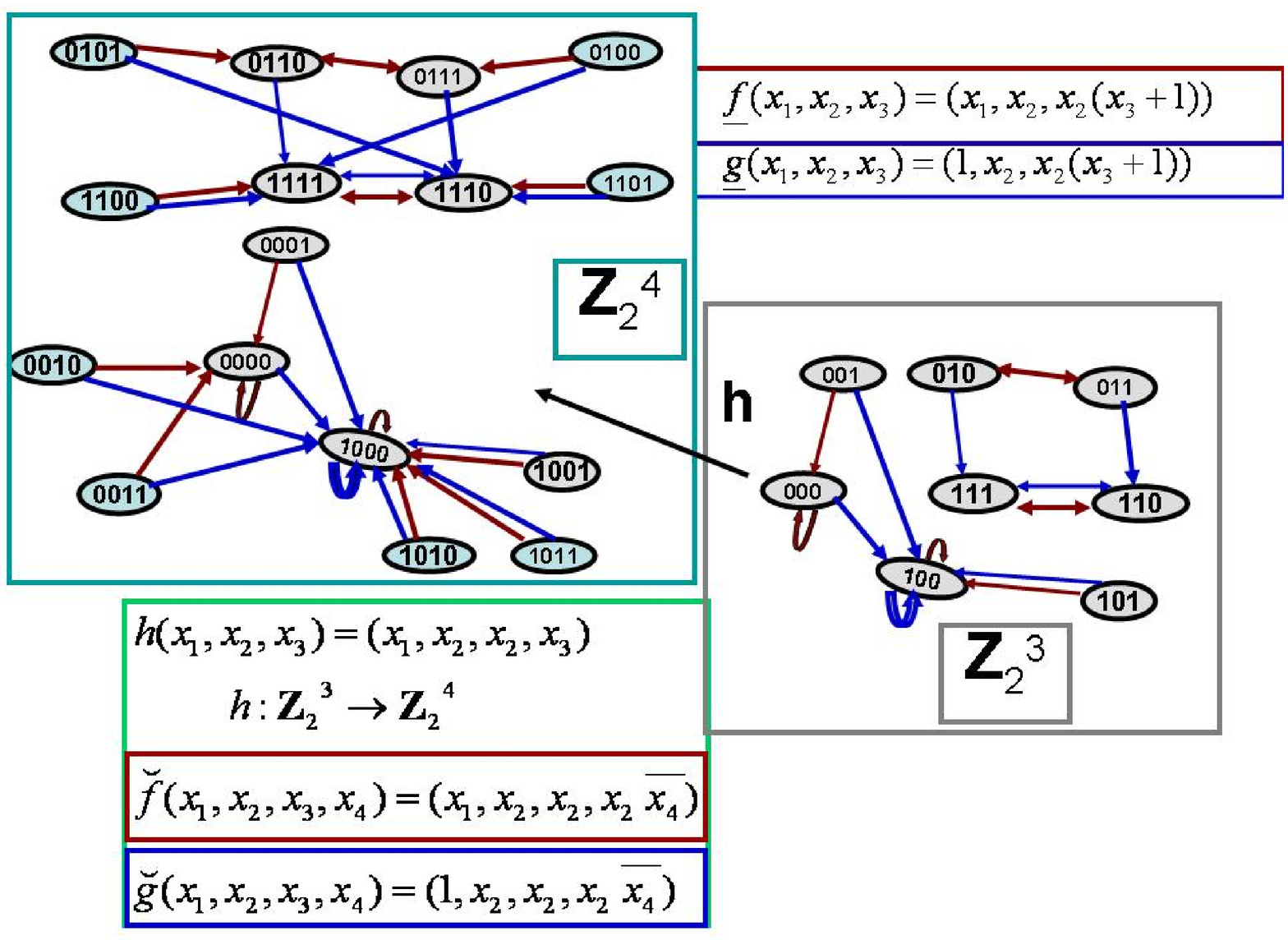}\caption{State Spaces of $\mathcal{F}$ and $\mathcal{F}'$. Transition Matrices: $T_1$, and $T_{h_\phi}$.}
\end{figure}
We verify the composition of  functions  as follows
\[\begin{array}{ccccccc}
 {\Z_p}^3 &\overset{f_{31}}{\longrightarrow} & {\Z_p}^3 &\overset{f_{21}}{\longrightarrow}& {\Z_p}^3
&\overset{f_{12}}{\longrightarrow} & {\Z_2}^3  \\
  h_\phi\downarrow & & h_\phi \downarrow& & h_\phi\downarrow &
&h_\phi\downarrow \\
 {\Z_p}^4 &\overset{g_{4}}{\longrightarrow} & {\Z_p}^4 &\overset{g_{21} \circ g_{32}}{\longrightarrow}
& {\Z_p}^4 &\overset{g_{11}}{\longrightarrow} & {\Z_p}^3
\end{array}\]
$(h_\phi \circ
f_{31})(x_1,x_2,x_3)=(x_1,x_2,x_2,x_2\overline{x_3})=(g_4\circ
h_\phi)(x_1,x_2,x_3)$, \\
$(h_\phi \circ f_{21})(x_1,x_2,x_3)=(x_1,x_2,x_2,x_3)=((g_{21} \circ
g_{32}) \circ
h_\phi)(x_1,x_2,x_3)$,\\
$(h_\phi \circ f_{12})(x_1,x_2,x_3)=(1,x_2,x_2,x_3)=(g_{11}\circ
h_\phi)(x_1,x_2,x_3)$\\
Similarly we check the condition for $\underline{f},$ and
$\check{f}$.

 The third condition holds, because the initial $\epsilon
 \leq .009$, because $|T_1-T_{h_\phi}|<.009$. In fact:
 \[T_1=\left( \begin{array}{llllllll}
.5168&0&0&0&.4832&0&0&0\\
.5168&0&0&0&.4832&0&0&0\\
0&0&0&.5168&0&0&0&.4832\\
0&0&.5168&0&0&0&.4832&0\\
0&0&0&0&1&0&0&0\\
 0&0&0&0&1&0&0&0\\
 0&0&0&0&0&0&0&1\\
  0&0&0&0&0&0&1&0\cr
\end{array}\right).\]
and
 \[T_{h_\phi}=\left( \begin{array}{llllllll}
.5168&0&0&0&.4832&0&0&0\\
.5168&0&0&0&.4832&0&0&0\\
0&0&0&.51428&0&0&0&.39999\\
0&0&0&.51428&0&0&0&.39999\\
.00252&0&0&0&.99748&0&0&0\\
 .00252&0&0&0&.99748&0&0&0\\
0&0&0&0&0&0&0&1\\
  0&0&0&0&0&0&1&0\cr
\end{array}\right).\]
To describe the connection between the  state spaces, see Figure 2.

\section{The category  \textbf{PSS}}
In this section, we  prove that the PSS with the
\emph{homomorphisms} form a category, that we denote by
\textbf{PSS}. In Theorem \ref{FSDS}, we prove that the category of
Sequential Dynamical Systems \textbf{SDS} is a full subcategory of
\textbf{PSS}.
\begin{definition2}\label{COM} Let
\[ \mathcal{H}_1=(\phi _1, h_{\phi _1}):\mathcal{F}\rightarrow \mathcal{G}\hbox{
and  } \mathcal{H}_2=(\phi _2, h_{\phi _2}):\mathcal{G}\rightarrow
\mathcal{L}\] be two PSS-homomorphisms. Then the composition
$\mathcal{H}=\mathcal{H}_2 \circ \mathcal{H}_1:
\mathcal{F}\rightarrow \mathcal{L}$  is defined as follows:
$\mathcal{H}=(\phi_1 \circ \phi _2, h_{\phi _2} \circ h_{\phi _1})$.

The map $\mathcal{H}: \mathcal{F}\rightarrow \mathcal{L}$ is a
PSS-homomorphism.
\end{definition2}
\begin{proof}
The composite $\phi=\phi_2 \circ \phi_1:\Lambda\rightarrow\Gamma $
of two graph morphisms is obviously again a graph morphism. The
composite $h_\phi=h_{\phi _2} \circ h_{\phi _1}$ is again a digraph
morphism which satisfies the conditions (\ref{homPSS}.1), and
(\ref{homPSS}.2). In fact, using the Proposition and Definition 2.7
in \cite{LP2}, we can check these condition and  conclude that the
third condition holds, too by Theorem \ref{MA}. So,  $(\phi,
h_\phi)$ is again a PSS-homomorphism .
\end{proof}
\begin{theorem}
The Probability Sequential Systems together with the homomorphisms
of PSS form the category \textbf{PSS}.
\end{theorem}
\begin{proof} Easily  follows from Definition and Theorem \ref{COM}.
\end{proof}
\begin{theorem}\label{FSDS} The SDS together with the morphisms defined in \ref{homPSS} form a full subcategory of the category
\textbf{PSS}.
\end{theorem}
\begin{proof}
It is trivial.
\end{proof}
In paper \cite{LP2}, the authors proved that the category SDS has
finite product, similarly the category \textbf{PSS} has finite
product too, so the following theorem for two PSS holds \cite{ML}

\begin{theorem}
Let $\mathcal{D}_S$, and $\mathcal{G}_{X}$ be two PSS over the
finite field $\mathcal{K}$. For all $\epsilon$-homomorphisms
$(\phi_1, h_{\phi_1}):\mathcal{L}\rightarrow\mathcal{D}_S$ and
$(\phi_2,h_{\phi_2}):\mathcal{L} \rightarrow \mathcal{G}_X$, then
there exists a morphism $(\phi, h_{\phi}):\mathcal{L} \rightarrow
\mathcal{D}_S \times \mathcal{G}_X$ such that the following diagram
commutes
\[\begin{array}{c}
  \hspace{.05in}\mathcal{D}_S \times
\mathcal{G}_X \hspace{.05in}\\
 \overset{\pi_1}{ \swarrow}\hspace{.1in}\overset{(\phi,h_{\phi})}{\uparrow}\hspace{.1in}\overset{\pi_2}{\searrow}\\
 \mathcal{D}_S\overset{(\phi_1,h_{\phi_1})}{\longleftarrow} \mathcal{L} \overset{(\phi_2,h_{\phi_2})}{\longrightarrow} \mathcal{G}_{X}\cr
  \end{array}\]
\end{theorem}

\subsection{The functor $ T: \textbf{PSS}\rightarrow \mathcal{A}b$ } Let
$\mathcal{D}=(\Gamma,\ (F_i)_{i=1}^{|\Gamma |=n},\ (\alpha _j)_j,\
S, C)$ be a PSS over the finite field $\mathcal{K}$. Let
$F=\{f:{\mathcal K}^{n}\rightarrow  {\mathcal
K}^{n}\}=S\dot{\cup}\overline{S}$ be the set of all  update
functions that we can construct with the local functions and the
permutation in $\mathcal{F}_S$. Let us consider, the free abelian
group generated by $F$, that we denote by $\langle F\rangle$, then
$\langle F \rangle=\langle S \rangle\oplus \langle \overline S
\rangle$. We can notice that we are working over a finite field with
characteristic the prime number $p$, so $pf=0$ for all $f\in S$. We
will take the quotients of these groups by $p$, that is  $\langle F
\rangle /\langle pF \rangle=\langle S \rangle/\langle pS
\rangle\oplus \langle \overline S \rangle /\langle p\overline{S}
\rangle$, and these groups are finite, \cite{KM,F,G}. Denoting
$\langle A \rangle /\langle pA \rangle=\langle A \rangle_p$, for an
abelian group $A$, we rewrite the above relation by  $\langle F
\rangle_p=\langle S \rangle_p\oplus \langle \overline S \rangle_p$
and, there exists a covariant functor from the category \textbf{PSS}
to the category of small abelian groups with morphism of such,
$\mathcal{A}b$, defined as follows.
\[ T: \textbf{PSS}\rightarrow \mathcal{A}b\]
1. The object function is defined by $ T (\mathcal{F}_S)=\langle S \rangle_p$.\\
2. The arrow function which assigns to each homomorphism $(\phi,
h_\phi): (\mathcal{F}_S)\rightarrow (\mathcal{G}_{X}) $ in the
category \textbf{PSS} an homomorphism of abelian groups  $ T
(\phi,h_\phi)=H_\phi:\langle S \rangle_p\rightarrow \langle X
\rangle_p$ which is defined in a natural way, because  $h_\phi\circ
(\sum_{f\in S}a f)=(\sum_{g\in S_\phi}a g)\circ h_\phi$, where $a\in
\Z_p$, and $S_\phi\subseteq X$, then
\[H_\phi(\sum_{f\in S}a f)=(\sum_{g\in X}a g).\]
$T$ is a functor , in fact,
\begin{itemize}
\item [1] $T(1_{\mathcal{F}_S})=1_{\langle S \rangle_p}$.
\item [2] $T((\phi_2,h_{\phi_2})\circ (\phi_1,h_{\phi_1}))=T((\phi_2,h_{\phi_2}))\circ T(
(\phi_1,h_{\phi_1}))$
\end{itemize}
The functor $T$ gives the possibility to work with PSS using  group
theory, for example, because $\langle F \rangle_p=\langle S
\rangle_p\oplus \langle \overline S\rangle _p$, and
\[\frac{\langle F \rangle _p}{\langle \overline S \rangle _p}\cong \langle  S \rangle _p,\]
We assign probabilities to the set $\overline S$ in some way, and we
consider that all  possible different assignations are
$\epsilon$-isomorphic.
\begin{definition}[(Complement of a PSS)] A complement of the PSS
$\mathcal{D}_S$ is the PSS  $\mathcal{D}_{\overline{S}}$, and all of
the complement are $\epsilon$-equivalent, so  we can select a
distribution of probabilities for the complement having in account
particular applications.
\end{definition}

We use the definition of complement in order to define a
decomposition of a PSS in two sub PSS, only looking the  set of
functions. One of the mean problem in modeling dynamical systems is
the computational aspect of the number of functions and the
calculation of steady states in the state space, in particular the
reduction of number of functions is one of the most important
problem to solve for determine which part of the network \emph{state
space} could be simplify.

\end{document}